\documentclass[final]{siamart171218}

\usepackage{amssymb}

\usepackage{algorithmic}
\usepackage[caption=false,font=footnotesize]{subfig}
\usepackage{color}
\usepackage{placeins}
\usepackage{tikz}

\usetikzlibrary{external}
\usepackage{pifont}
\newcommand{\cmark}{\ding{51}}%
\newcommand{\xmark}{\ding{55}}%

\newcommand{\mat}[1]            {\boldsymbol{#1}}

\newcommand{\norm}[1]           {\left|\!\left|#1\right|\!\right|}
\renewcommand{\tilde}[1]        {\widetilde{#1}}
\renewcommand{\hat}[1]          {\widehat{#1}}
\newcommand{\steering}{\sigma}
\newcommand{\field}[1]{\mathbb{#1}}
\newcommand{\N}{\field{N}}

\newcommand{\T}{\mathcal{T}}
\newcommand{\Xt}{\widetilde{X}}
\newcommand{\Tt}{\widetilde{\T}}

\newcommand{\True}{\texttt{True}}
\newcommand{\False}{\texttt{False}}

\renewcommand{\vec}[1]{#1}

\usepackage{eqparbox}

\usepackage{url}
\DeclareUrlCommand\UScore{\urlstyle{tt}}
\newcommand{\mpi}[1]{\UScore{MPI_#1}}

\usepackage{cleveref}
\crefname{algorithm}{algorithm}{algorithms}
\Crefname{algorithm}{Algorithm}{Algorithms}
\crefname{subsection}{section}{sections}
\Crefname{subsection}{Section}{Sections}

\begin{document}

\title{Scalable Asynchronous\\ Domain Decomposition Solvers\thanks{
    Part of this work has been accepted for publication in the form of a proceedings paper by the 25th International Domain Decomposition Conference.
  }}

\author{
  Christian Glusa\thanks{Center for Computing Research, Sandia National Laboratories, Albuquerque, New Mexico, USA ({\tt caglusa@sandia.gov}).}
  \and
  Erik G. Boman\thanks{Center for Computing Research, Sandia National Laboratories, Albuquerque, New Mexico, USA ({\tt egboman@sandia.gov}).}
  \and
  Edmond Chow\thanks{School of Computational Science and Engineering, College of Computing, Georgia Institute of Technology, Atlanta, Georgia, USA ({\texttt{echow@cc.gatech.edu}}).}
  \and
  Sivasankaran Rajamanickam\thanks{Center for Computing Research, Sandia National Laboratories, Albuquerque, New Mexico, USA ({\tt srajama@sandia.gov}).}
  \and
  Daniel B. Szyld\thanks{Temple University, Philadelphia, Pennsylvania, USA ({\texttt{szyld@temple.edu}}).}
}

\maketitle

\begin{abstract}
  Parallel implementations of linear iterative solvers generally alternate between phases of data exchange and phases of local computation.
  Increasingly large problem sizes and more heterogeneous compute architectures make load balancing and the design of low latency network interconnects that are able to satisfy the communication requirements of linear solvers very challenging tasks.
  In particular, global communication patterns such as inner products become increasingly limiting at scale.

  We explore the use of asynchronous communication based on one-sided MPI primitives in the context of domain decomposition solvers.
  In particular, a scalable asynchronous two-level Schwarz method is presented.
  We discuss practical issues encountered in the development of a scalable solver and show experimental results obtained on a state-of-the-art supercomputer system that illustrate the benefits of asynchronous solvers in load balanced as well as load imbalanced scenarios.
  Using the novel method, we can observe speed-ups of up to 4x over its classical synchronous equivalent. 
\end{abstract}



\section{Introduction}

Multilevel methods such as multigrid and domain decomposition are among the most efficient and scalable solvers for partial differential equations developed to date.
Adapting them to the next generation of supercomputers and improving their performance and scalability is crucial in the push towards exascale.
Domain decomposition methods subdivide the global problem into subdomains, and then alternate between local solves and boundary data exchange.
This puts a significant stress on the network interconnect, since all processes try to communicate at once.
On the other hand, during the solve phase, the network is under-utilized.
The use of non-blocking communication can only alleviate this issue, but not fully resolve it.
In asynchronous methods, on the other hand, computation and communication occur at the same time, with some processes performing computation while others communicate, so that the network is consistently in use.

The term ``asynchronous'' can have several different meanings in the literature.
In computer science, it is sometimes used to describe communication patterns that are non-blocking, so that computation and communication can be overlapped.
Iterative algorithms that use such ``asynchronous'' communication yield the same iterates (results) up to round-off error, as they do not change the mathematical algorithm.
In applied mathematics, on the other hand, ``asynchronous'' denotes parallel algorithms where each process (processor) proceeds at its own speed without synchronization.
Thus, asynchronous algorithms go beyond the widely used bulk-synchronous parallel (BSP) model.
More importantly, they are mathematically different than synchronous methods and generate different iterates.
The earliest work in this area was called ``chaotic relaxation'' \cite{ChazanMiranker1969_ChaoticRelaxation}.
Both types of asynchronous approaches are expected to play an important role on future supercomputers.
In this paper, we focus on asynchronous methods in the mathematical sense, and we will use the terms ``asynchronous'' and ``synchronous'' to distinguish between methods that are asynchronous and synchronous in the mathematical sense.

Domain decomposition solvers \cite{DoleanJolivetEtAl2015_IntroductionToDomainDecompositionMethods,ToselliWidlund2006_DomainDecompositionMethods,SmithBjorstadEtAl2004_DomainDecomposition} are often used as preconditioners in Krylov subspace iterations.
Unfortunately, the computation of inner products and norms widely used in Krylov methods requires global communication.
Global communication primitives, such as \mpi{Reduce}, asymptotically scale as the logarithm of the number of processes involved.
This can become a limiting factor when very large process counts are used.
The underlying domain decomposition method, however, can do away with globally synchronous communication, assuming the coarse problem in multilevel methods can be solved in a parallel way.
Therefore, we will focus on using domain decomposition methods purely as iterative methods in the present work.
We will note, however, that the discussed algorithms could be coupled with existing pipelined methods \cite{GhyselsAshbyEtAl2013_HidingGlobalCommunicationLatency} which alleviate the global synchronization requirement of Krylov solvers.

Another issue that is crucial to good scaling behavior is load imbalance.
Load imbalance might occur due to heterogeneous hardware in the system, network noise, dynamic power capping \cite{AhlgrenAnderssonEtAl2018_CraySystemMonitoring}, or due to local, problem specific causes, such as iteration counts for local solves that vary from subdomain to subdomain.
The latter are especially difficult to predict, so that load balancing cannot occur before the actual solve.
Therefore, processes in a synchronous parallel program must be idle until its slowest process has finished.
In an asynchronous method, local computation can continue, and potentially improve the quality of the global solution.

An added benefit of asynchronous methods is that, since the interdependence between subdomains has been weakened, fault tolerance \cite{CappelloGeistEtAl2009_TowardExascaleResilience,CappelloGeistEtAl2014_TowardExascaleResilience} can be more easily achieved.
When one process must stop, be it for a hard or a soft fault, it can be replaced without having to halt every other process.

The main drawback of asynchronous iterations is the fact that deterministic behavior is sacrificed.
Consecutive runs do not produce the same result.
(But one would hope that they are at most a distance proportional to the convergence tolerance apart from each other.)
This also makes the mathematical analysis of asynchronous methods significantly more difficult than for its synchronous counterparts.
Analytical frameworks for asynchronous linear (and nonlinear) iterations have long been available \cite{ChazanMiranker1969_ChaoticRelaxation,Baudet1978_AsynchronousIterativeMethodsMultiprocessors,Bertsekas1983_DistributedAsynchronousComputationFixedPoints,FrommerSzyld2000_AsynchronousIterations}, but generally cannot produce sharp convergence bounds except in the simplest cases.

The main contributions of our work are:
\begin{itemize}
\item A novel asynchronous two-level domain decomposition method, scalable to thousands of processors.
\item An empirical study of one-sided MPI performance in a scientific computing setting.
\item Empirical comparisons of synchronous and asynchronous variants of domain decomposition solvers on a state-of-the-art parallel computer.
\end{itemize}

Our work demonstrates that asynchronous methods have the potential of outperforming conventional synchronous solvers and offer a viable alternative in the push towards exascale.

The present work is structured as follows:
In \Cref{sec:dd-methods}, we present overlapping domain decomposition methods, and explain their use in synchronous and asynchronous fashion.
For a general introduction to domain decomposition methods we refer the reader to \cite{DoleanJolivetEtAl2015_IntroductionToDomainDecompositionMethods,ToselliWidlund2006_DomainDecompositionMethods,SmithBjorstadEtAl2004_DomainDecomposition}.
The section concludes with a convergence analysis of the presented one- and two-level methods.
\Cref{sec:one-sided-mpi} is dedicated to a description of the presently available mechanisms in MPI and hardware to achieve truly asynchronous communication.
Numerical experiments exploring asynchronous communication and using the presented domain decomposition methods are given in \Cref{sec:experiments}, where we compare the strong and weak scaling behavior of synchronous and asynchronous solvers with and without load imbalance.

\subsection{Related work}

An asynchronous one-level domain decomposition solver with optimized artificial boundary conditions was proposed in \cite{MagoulesSzyldEtAl2017_AsynchronousOptimizedSchwarzMethods}; see also \cite{GarayMagoulesEtAl2018_ConvergenceAsynchronousOptimizedSchwarzMethodsPlane,GarayMagoulesEtAl2017_SynchronousAsynchronousOptimizedSchwarz,HaddadGarayEtAl2018_SynchronousAsynchronousOptimizedSchwarz} for its analysis in two different settings.
An implementation of asynchronous optimized Schwarz is described in \cite{YamazakiChowEtAl2019_PerformanceAsynchronousOptimizedSchwarz}.
An optimization package that leverages asynchronous coordinate updates is presented in \cite{PengXuEtAl2016_ARock}.
An asynchronous multigrid method for shared memory systems was proposed in \cite{Wolfson-PouChow2019_AsynchronousMultigridMethods}.
Synchronization reducing Krylov methods have a long history~\cite{chronopoulos1989step}.
However, preconditioning such methods is unresolved apart from some simple preconditioners~\cite{chronopoulos1989implementation}.
Recent work extends their applicability to one level domain decomposition preconditioning~\cite{yamazaki2014domain}.
Pipelined Krylov methods~\cite{GhyselsAshbyEtAl2013_HidingGlobalCommunicationLatency} reduce synchronization costs by overlapping inner products with matrix-vector products and preconditioner applications, and can be used with any preconditioner.


\section{Domain decomposition methods}
\label{sec:dd-methods}

\subsection{One-level Restricted Additive Schwarz (RAS)}

We want to solve the global system
\begin{align*}
  \mat{A}\vec{u}=\vec{f},
\end{align*}
where \(\mat{A}\in\mathbb{R}^{N\times N}\) arises from the finite element or finite difference discretization of a partial differential equation.
Informally, one-level domain decomposition solvers break up the global system of equations into overlapping sub-problems that cover the whole global system.
This requires that the matrix \(\mat{A}\) is sparse and couples unknowns only in a local manner.

The iteration then alternates between computation of the global residual, which involves communication, and local solves for solution corrections.
Special attention needs to be paid to the unknowns in the overlap, in order to avoid over-correction.
Below, we describe the different methods considered in this work in detail in order to understand what data is required to be exchanged and how the methods can be executed in asynchronous fashion.

Based on the graph of $\mat{A}$ or geometric information for the underlying problem the unknowns are grouped into \(P\) overlapping sets \(\mathcal{N}_{p}\) of size \(N_{p}\), \(p=1, \ldots, P\).
An example of such a partitioning is given in \Cref{fig:partitioning}.
We further split the sets \(\mathcal{N}_{p}\) into 
\begin{align*}
  \mathcal{S}_{p} := \left\{j\in\mathcal{N}_{p} \mid \exists k\in\mathcal{N}_{p}^{c}: \mat{A}_{jk}\neq 0\right\},
\end{align*}
i.e., unknowns that are on the boundary of the set \(\mathcal{N}_{p}\), and interior unknowns \(\mathcal{I}_{p}:=\mathcal{N}_{p}\setminus\mathcal{S}_{p}\).

{
  \setlength{\abovecaptionskip}{-10pt}
  \begin{figure}[!t]
    \centering
    \includegraphics[width=0.5\columnwidth]{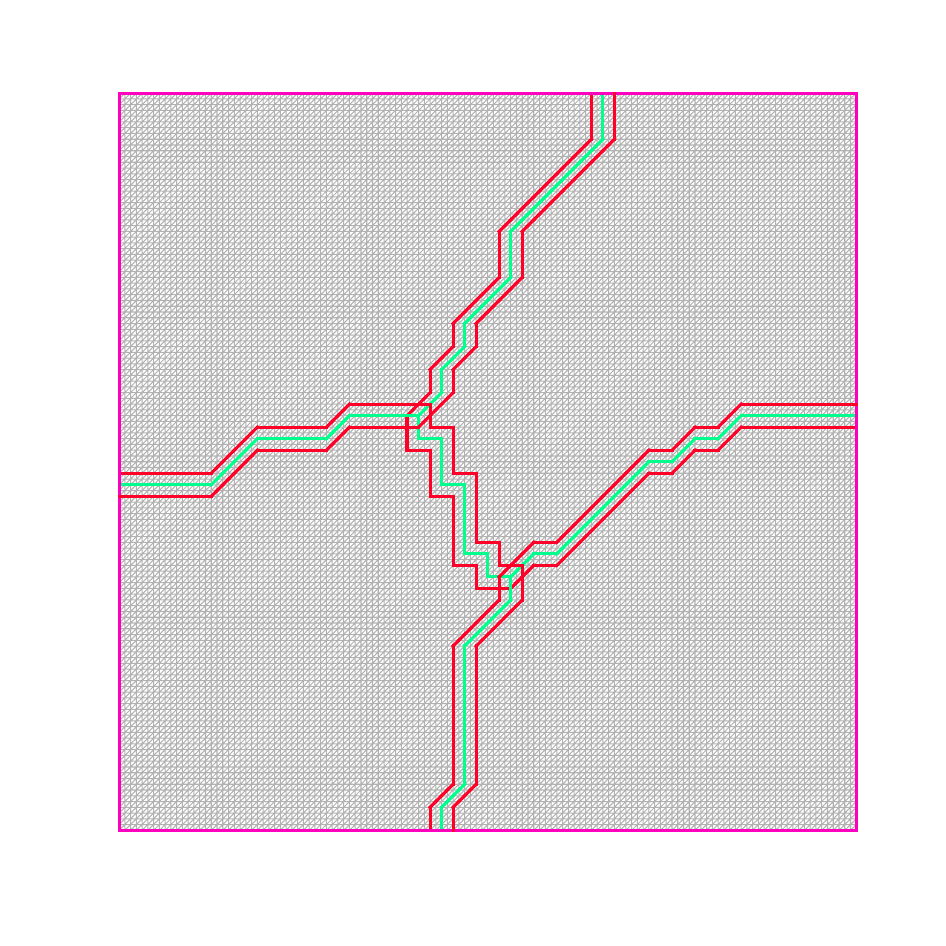}
    \caption{
      Partitioning of a uniform triangular mesh of the unit square into 4 overlapping subdomains.
      The non-overlapping partitioning produced using METIS \cite{KarypisKumar1998_FastHighQualityMultilevel} is shown in green; the extended overlapping subdomains are shown in red.
    }\label{fig:partitioning}
  \end{figure}
}
The notation throughout this section is based on Dolean et al.~\cite{DoleanJolivetEtAl2015_IntroductionToDomainDecompositionMethods}.
We call the restriction to the \(p\)-th set \(\mat{R}_{p}\in\mathbb{R}^{N_{p}\times N}\).
The entries of the matrices \(\mat{R}_{p}\) are all either one or zero, with exactly one entry per row and at most one entry per column being non-zero.
The local parts of \(\mat{A}\) are given by
\begin{align*}
  \mat{A}_{p}&= \mat{R}_{p}\mat{A}\mat{R}_{p}^{T} \in\mathbb{R}^{N_{p}\times N_{p}}.
\end{align*}

Furthermore, we require a partition of unity, represented by diagonal weighting matrices \(\mat{D}_{p}\), such that the discrete partition of unity property holds
\begin{align}
  \mat{I}=\sum_{p=1}^{P}\mat{R}_{p}^{T}\mat{D}_{p}\mat{R}_{p}.\label{eq:partitionOfUnity}
\end{align}
In what follows, we will assume that \(\mat{D}_{p}\) are Boolean, i.e. their entries are either zero or one.
This means that every (potentially shared) unknown has a special attachment with exactly one subdomain.
We will furthermore require that \(\left(\mat{D}_{p}\right)_{jj}=0\) for all surface unknowns \(j\in\mathcal{S}_{p}\).
One way of satisfying these restrictions is to extend overlaps starting with a \emph{non-overlapping} partition and then define the special attachment via the partition.

\begin{figure}
  \centering
  \includegraphics{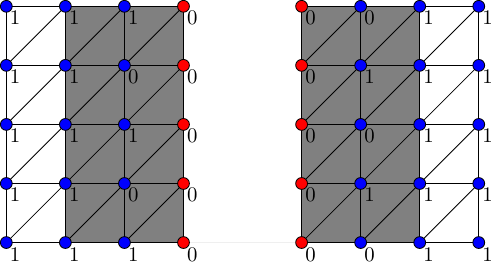}
  \caption{
    Two overlapping subdomains.
    The overlap between the subdomains is shaded in gray; the respective surface sets \(\mathcal{S}_{\bullet}\) are shown by red circles, the interior unknowns \(\mathcal{I}_{\bullet}\) as blue circle.
    The diagonal values of the respective \(\mat{D}_{\bullet}\) are shown next to the nodes.
  }
\end{figure}

Consequently,
\begin{align}
  \mat{D}_{p}\mat{R}_{p}\mat{R}_{q}^{T}\mat{D}_{q}&=\mat{0} \quad \text{for } p\neq q \label{eq:2}
\end{align}
and
\begin{align}
  \label{eq:3}
  \mat{D}_{p}\mat{R}_{p}\mat{R}_{p}^{T}\mat{D}_{p} &= \mat{D}_{p}.
\end{align}
Moreover, the identity 
\begin{align}
  \mat{R}_{p}\mat{A}\mat{R}_{q}^{T}\mat{D}_{q}
  = \mat{R}_{p}\mat{R}_{q}^{T}\mat{R}_{q}\mat{A}\mat{R}_{q}^{T}\mat{D}_{q} \label{eq:1}
\end{align}
holds, since for any \(\vec{u}_{q}\in\mathbb{R}^{N_{q}}\), \(\mat{D}_{q}\vec{u}_{q}\) is supported on the interior unknowns \(\mathcal{I}_{p}\), and hence \(\mat{A}\mat{R}_{q}^{T}\mat{D}_{q}\vec{u}_{q}\) is supported in \(\mathcal{N}_{q}\).
But on \(\mathcal{N}_{q}\), \(\mat{R}_{q}^{T}\mat{R}_{q}\) acts as the identity.

A stationary iterative method based on the
splitting $\mat{A} = \mat{M} - \mat{N}$ is given globally as
\begin{align*}
  \vec{u}^{n+1}&=\vec{u}^{n} + \mat{M}^{-1} \left(\vec{f}-\mat{A}\vec{u}^{n}\right),
\end{align*}
where \(\mat{M}^{-1}\) is a preconditioner for \(\mat{A}\).

This means that we need to calculate the residual \(\vec{r}^{n}=\vec{f}-\mat{A}\vec{u}^{n}\).
Its local part on node \(p\) is given by
\begin{align*}
  \mat{R}_{p}\vec{r}^{n}
  &= \mat{R}_{p}\vec{f} - \mat{R}_{p}\mat{A}\vec{u}^{n} \\
  &= \mat{R}_{p}\left(\sum_{q=1}^{P}\mat{R}_{q}^{T}\mat{D}_{q}\mat{R}_{q}\right)\vec{f} - \mat{R}_{p}\mat{A}\left(\sum_{q=1}^{P}\mat{R}_{q}^{T}\mat{D}_{q}\mat{R}_{q}\right)\vec{u}^{n} \\
  &= \sum_{q=1}^{P}\mat{R}_{p}\mat{R}_{q}^{T}\mat{D}_{q}\mat{R}_{q}\vec{f} - \sum_{q=1}^{P}\mat{R}_{p}\mat{R}_{q}^{T}\mat{A}_{q}\mat{D}_{q}\mat{R}_{q}\vec{u}^{n}\\
  &=\sum_{q=1}^{P}\mat{R}_{p}\mat{R}_{q}^{T}\left(\mat{D}_{q}\mat{R}_{q}\vec{f} - \mat{A}_{q}\mat{D}_{q}\mat{R}_{q}\vec{u}^{n}\right),
\end{align*}
where we used~\eqref{eq:partitionOfUnity} and \eqref{eq:1}.
This means that in order to obtain the local part of the global residual, we first compute locally \(\mat{D}_{p}\mat{R}_{p}\vec{f} - \mat{A}_{p}\mat{D}_{p}\mat{R}_{p}\vec{u}^{n}\) on every node \(p\), and then communicate and accumulate the overlapping parts of these local residual vectors.
The latter operation is represented by the operator \(\sum_{q=1}^{P}\mat{R}_{p}\mat{R}_{q}^{T}\).

The \emph{restricted additive Schwarz (RAS) preconditioner} \cite{CaiSarkis1999_RestrictedAdditiveSchwarzPreconditioner,CaiDryjaEtAl2003_RestrictedAdditiveSchwarzPreconditioners} is given by
\begin{align*}
  \mat{M}_{RAS}^{-1} = \sum_{p=1}^{P}\mat{R}_{p}^{T}\mat{D}_{p}\mat{A}_{p}^{-1}\mat{R}_{p}.
\end{align*}
RAS is widely used and is the default option for overlapping domain decomposition preconditioners in PETSc \cite{petsc-user-ref}.
It can be thought of as a variant of the additive Schwarz preconditioner
\begin{align*}
  \mat{M}_{AS}^{-1} = \sum_{p=1}^{P}\mat{R}_{p}^{T}\mat{A}_{p}^{-1}\mat{R}_{p}
\end{align*}
that is convergent as an iterative method, since the damping by \(\mat{D}_{p}\) in the overlapping parts avoids over-correction; see \cite{FrommerSzyld2000_AsynchronousIterations}.
Note that for a natural choice of $\mat{D}_{p}$, the number of communication steps is cut in half as there is no communication associated with $\mat{R}^T_p \mat{D}_p$.

Now, the local part of the RAS iteration is given by
\begin{align*}
  \mat{R}_{p}\vec{u}^{n+1}
  &= \mat{R}_{p}\vec{u}^{n} + \mat{R}_{p}\mat{M}_{RAS}^{-1}\vec{r}^{n}\\
  &= \mat{R}_{p}\vec{u}^{n} + \sum_{q=1}^{P} \mat{R}_{p}\mat{R}_{q}^{T}\mat{D}_{q}\mat{A}_{q}^{-1}\mat{R}_{q}\vec{r}^{n}.
\end{align*}
If we set \(\vec{u}_{p}^{n}=\mat{R}_{p}\vec{u}^{n}\) and \(\vec{r}_{p}^{n}=\mat{R}_{p}\vec{r}^{n}\) as the local parts of solution and residual respectively, the RAS iteration is
\begin{align*}
  \vec{r}_{p}^{n} &= \sum_{q=1}^{P} \mat{R}_{p}\mat{R}_{q}^{T}\left(\mat{D}_{q}\mat{R}_{q}\vec{f} - \mat{A}_{q}\mat{D}_{q}\vec{u}_{q}^{n}\right), \\
  \vec{u}_{p}^{n+1} &= \vec{u}_{p}^{n} + \sum_{q=1}^{P} \mat{R}_{p}\mat{R}_{q}^{T}\mat{D}_{q}\mat{A}_{q}^{-1}\vec{r}_{q}^{n}.
\end{align*}
This seems to suggest that the update step requires neighborhood communication as well.
But in fact, in the next iteration, computation of the residual only requires \(\mat{D}_{p}\vec{u}_{p}^{n+1}\).
From \eqref{eq:2}, \eqref{eq:3}, we see that the iterative scheme without the communication step in the update
\begin{align}
  \vec{r}_{p}^{n} &= \sum_{q=1}^{P} \mat{R}_{p}\mat{R}_{q}^{T}\left(\mat{D}_{q}\mat{R}_{q}\vec{f} - \mat{A}_{q}\mat{D}_{q}\vec{w}_{q}^{n}\right), \label{eq:4}\\
  \vec{w}_{p}^{n+1} &= \vec{w}_{p}^{n} + \mat{A}_{p}^{-1}\vec{r}_{p}^{n}\label{eq:5}
\end{align}
is equivalent because \(\mat{D}_{p}\vec{u}_{p}^{n}=\mat{D}_{p}\vec{w}_{p}^{n}\) for all \(n\).
The solution \(\vec{u}_{p}^{n}\) can be recovered from \(\vec{w}_{p}^{n}\) in the post-processing step
\begin{align*}
  \vec{u}_{p}^{n} &= \mat{R}_{p}\vec{u}^{n}=\sum_{q=1}^{P}\mat{R}_{p}\mat{R}_{q}^{T}\mat{D}_{q}\mat{R}_{q}\vec{u}^{n} = \sum_{q=1}^{P}\mat{R}_{p}\mat{R}_{q}^{T}\mat{D}_{q}\vec{w}_{q}^{n}.
\end{align*}

Finally, we use the norm of the residual in the stopping criterion.
The norm can be computed from local quantities as
\begin{align*}
  \norm{\vec{r}^{n}}^{2}
  &= \vec{r}^{n} \cdot \vec{r}^{n} = \vec{r}^{n} \cdot \left(\sum_{p=1}^{P}\mat{R}_{p}^{T}\mat{D}_{p}\mat{R}_{p}\vec{r}^{n}\right)\\
  &=\sum_{p=1}^{P}\left(\mat{R}_{p}\vec{r}^{n}\right) \cdot \left(\mat{D}_{p}\mat{R}_{p}\vec{r}^{n}\right) 
  =\sum_{p=1}^{P}\vec{r}_{p}^{n} \cdot \left(\mat{D}_{p}\vec{r}_{p}^{n}\right).
\end{align*}

In conclusion, we can give the local form of RAS as in \Cref{alg:RAS}, where we have dropped the superscript \(n\) for the iteration number.
In fact, \Cref{alg:RAS} describes both the synchronous \emph{and} the asynchronous version of RAS.
In the synchronous version, line 4 is executed in lock step fashion by all subdomains using non-blocking two-sided communication primitives.
This communication step could be overlapped by computation.
However, in established frameworks such as Trilinos, such overlapping requires major changes to the framework\footnote{https://github.com/trilinos/Trilinos/issues/767}.
PETSc allows some overlap of computation and communication with two-phase assembly~\cite{petsc-user-ref}.
It is possible to modify such established libraries for the asynchronous iterations of this paper.
However, in order to keep the focus on algorithmic development, we developed a library that supports the one-sided communication primitives, and build the new solvers using the communication primitives.

In the asynchronous variant, each subdomain exposes a memory region for remote access.
On execution of line 4, the relevant components of the current local residual vector \(\vec{s}_{p}=\mat{D}_{p}\mat{R}_{p}\vec{f}-\mat{A}_{p}\mat{D}_{p}\vec{w}_{p}\) are written to the neighboring subdomains, and the latest locally available data \(\vec{s}_{q}\) from every neighbor \(q\) is used.
We refer to \Cref{sec:selection-one-sided} for a discussion of the options for actually achieving this neighborhood exchange in practice.
The implementation of a convergence check (as used on line 2) that does not require synchronization is detailed in \Cref{sec:conv-detection}.

\begin{algorithm}[!t]
  \begin{algorithmic}[1]
    \STATE \(\vec{w}_{p}\leftarrow \vec{0}\)
    \WHILE{not converged}
    \STATE Local residual: \(\vec{s}_{p}\leftarrow\mat{D}_{p}\mat{R}_{p}\vec{f} - \mat{A}_{p}\mat{D}_{p}\vec{w}_{p}\)
    \STATE Accumulate: \(\vec{r}_{p}\leftarrow\sum_{q=1}^{P}\mat{R}_{p}\mat{R}_{q}^{T}\vec{s}_{q}\) \COMMENT{\(\leftrightarrow\)}
    \STATE Solve: \(\mat{A}_{p}\vec{v}_{p}=\vec{r}_{p}\)
    \STATE Update: \(\vec{w}_{p}\leftarrow \vec{w}_{p}+ \vec{v}_{p}\)
    \ENDWHILE
    \STATE Post-process: \(\vec{u}_{p}\leftarrow\sum_{q=1}^{P}\mat{R}_{p}\mat{R}_{q}^{T}\mat{D}_{q}\vec{w}_{q}\) \COMMENT{\(\leftrightarrow\)}
  \end{algorithmic}
  \caption{Restricted additive Schwarz (RAS) in local form, ``\(\leftrightarrow\)'' signifies communication.}\label{alg:RAS}
\end{algorithm}

\subsection{Two-level synchronous RAS}

In order to improve the scalability of the solver, a mechanism of global information exchange is required.
Let \(\mat{R}_{0}\in\mathbb{R}^{n_{0}\times n}\) be the restriction from the fine grid problem to a coarser mesh, and let the coarse-grid matrix \(\mat{A}_{0}\) be given by the Galerkin relation \(\mat{A}_{0}=\mat{R}_{0}\mat{A}\mat{R}_{0}^{T}\).
The coarse-grid solve can be incorporated in the RAS iteration either in additive fashion:
\begin{align}
  \vec{u}^{n+1}&=\vec{u}^{n} + \left(\frac{1}{2}\mat{M}_{RAS}^{-1} + \frac{1}{2}\mat{R}_{0}^{T}\mat{A}_{0}^{-1}\mat{R}_{0}\right) \left(\vec{f}-\mat{A}\vec{u}^{n}\right), \label{eq:additiveRAS}
\end{align}
or in multiplicative fashion:
\begin{align*}
  \vec{u}^{n+1/2}&=\vec{u}^{n} + \mat{R}_{0}^{T}\mat{A}_{0}^{-1}\mat{R}_{0} \left(\vec{f}-\mat{A}\vec{u}^{n}\right), \\
  \vec{u}^{n+1}&=\vec{u}^{n+1/2} + \mat{M}_{RAS}^{-1}\left(\vec{f}-\mat{A}\vec{u}^{n+1/2}\right).
\end{align*}
In what follows, we focus on the additive version, since it naturally lends itself to asynchronous iterations: subdomain solves and coarse-grid solves are independent of each other.


We now determine the local form of the global algorithm.
It is understood that the solve with \(\mat{A}_{0}\) itself might be distributed over several processes.
This internal computation is not meant to be performed in an asynchronous manner, which is why we do not need to further explore the local form of the coarse-grid solve.
For simplicity of exposition we therefore do not describe the solution of the coarse-grid problem itself in local form, i.e. we will simply write \(\mat{A}_{0}^{-1}\).
The local part of the coarse-grid update is
\begin{align*}
  &\frac{1}{2}\mat{R}_{p}\mat{R}_{0}^{T}\mat{A}_{0}^{-1}\mat{R}_{0}\left(\vec{f}-\mat{A}\vec{u}^{n}\right) \\
  =&\frac{1}{2}\left(\mat{R}_{p}\mat{R}_{0}^{T}\right)\mat{A}_{0}^{-1}\sum_{p=1}^{P}\left(\mat{R}_{0} \mat{R}_{p}^{T}\right)\left(\mat{D}_{p}\mat{R}_{p}\vec{f} - \mat{A}_{p}\mat{R}_{p}\vec{u^{n}}\right).
\end{align*}
Here, the operators \(\left(\mat{R}_{0}\mat{R}_{p}^{T}\right)\) and \(\left(\mat{R}_{p}\mat{R}_{0}^{T}\right)\) encode the communication from subdomain \(p\) to the coarse grid and vice versa.
We notice that while the communication among subdomains consist in one neighborhood data exchange per iteration, the coarse-grid solve involves sending data from the subdomains to the coarse grid, and sending a solution from the coarse grid to the subdomains.
In conclusion, the local form of RAS with an additive coarse grid is given in \Cref{alg:additiveRAS}.
Again, we have dropped the superscript for the iteration number.
The communication between coarse and fine grid can be implemented in multiple ways.
Since we want to allow the coarse grid solve to be distributed itself and the same coarse unknown can be owned by several coarse grid ranks (just as is the case for the fine grid), we do not consider options involving \mpi{Reduce}/\mpi{Bcast} or \mpi{Gather}/\mpi{Scatter} or their non-blocking equivalents.
Instead, we opted for use of \mpi{Isend} and \mpi{Irecv}.
A future improvement could involve the use of intercommunicators and \mpi{Iallgatherv} or other collectives.
The advantage of the current approach is that the changes between synchronous and asynchronous implementation of the communication layer (described in the next section) are minimal.

\begin{algorithm}[!t]
  \begin{algorithmic}[1]
    \newcommand{\Subdomain}{\STATE \textbf{On subdomains} \begin{ALC@g}}
    \newcommand{\EndSubdomain}{\end{ALC@g}}
    \newcommand{\Coarse}{\STATE \textbf{On coarse grid} \begin{ALC@g}}
    \newcommand{\EndCoarse}{\end{ALC@g}}
    \STATE \(\vec{w}_{p}\leftarrow \vec{0}\)
    \WHILE{not converged}
    \Subdomain
    \STATE Local residual: \(\vec{s}_{p}\leftarrow\mat{D}_{p}\mat{R}_{p}\vec{f} - \mat{A}_{p}\mat{D}_{p}\vec{w}_{p}\)
    \STATE Send \(\mat{R}_{0}\mat{R}_{p}^{T}\vec{s}_{p}\) to coarse grid \COMMENT{\(\leftrightarrow\)}
    \STATE Accumulate: \(\vec{r}_{p}\leftarrow\sum_{q=1}^{P}\mat{R}_{p}\mat{R}_{q}^{T}\vec{s}_{q}\) \COMMENT{\(\leftrightarrow\)}
    \STATE Solve: \(\mat{A}_{p}\vec{v}_{p}=\vec{r}_{p}\)
    \STATE Update: \(\vec{w}_{p}\leftarrow \vec{w}_{p}+ \frac{1}{2}\vec{v}_{p}\)
    \STATE Receive \(\vec{c}_{p}=\mat{R}_{p}\mat{R}_{0}^{T}\vec{v}_{0}\) from coarse grid \COMMENT{\(\leftrightarrow\)}
    \STATE Update: \(\vec{w}_{p}\leftarrow \vec{w}_{p}+ \frac{1}{2}\vec{c}_{p}\)
    \EndSubdomain
    \Coarse
    \STATE Receive \(\mat{R}_{0}\mat{R}_{p}^{T}\vec{s}_{p}\) from subdomains \COMMENT{\(\leftrightarrow\)}
    \STATE Accumulate \(\vec{r}_{0}=\sum_{p=1}^{P}\mat{R}_{0}\mat{R}_{p}^{T}\vec{s}_{p}\)
    \STATE Solve \(\mat{A}_{0}\vec{v}_{0}=\vec{r}_{0}\)
    \STATE Send \(\vec{c}_{p}=\mat{R}_{p}\mat{R}_{0}^{T}\vec{v}_{0}\), \(p=1,\dots,P\) to subdomains \COMMENT{\(\leftrightarrow\)}
    \EndCoarse
    \ENDWHILE
    \Subdomain
    \STATE Post-process \(\vec{u}_{p}\leftarrow\sum_{q=1}^{P}\mat{R}_{p}\mat{R}_{q}^{T}\mat{D}_{q}\vec{w}_{q}\) \COMMENT{\(\leftrightarrow\)}
    \EndSubdomain
\end{algorithmic}
\caption{Synchronous RAS with additive coarse grid in local form, ``\(\leftrightarrow\)'' signifies communication.}\label{alg:additiveRAS}
\end{algorithm}

\subsection{Two-level asynchronous RAS}

From the mathematical description \eqref{eq:additiveRAS} of two-level additive RAS, one might be tempted to see the coarse-grid problem simply as an additional subdomain.
From \Cref{alg:additiveRAS} the fundamental differences between the subdomains and the coarse-grid problem become apparent.
Subdomains determine the right-hand side for their local solve and correct it by transmitting boundary data to their neighbors.
The coarse grid, on the other hand, receives its entire right-hand side from the subdomains, and hence it has to communicate with every single one of them.

In order to perform asynchronous coarse-grid solves, we therefore need to make sure that all the right-hand side data necessary for the solve has been received by the processes responsible for the coarse grid.
Moreover, corrections sent by the coarse grid should be used exactly once by the subdomains.
This is achieved by not only allocating memory regions to hold the coarse-grid right-hand side on the coarse-grid processes and the coarse-grid correction on the subdomains, but also Boolean variables that are polled to determine whether writing or reading right-hand side or solution data is permitted.
More precisely, writing of the local subdomain residuals to the coarse-grid memory region of \({\color{red} \vec{r}_{0}}\) is contingent upon the state of the Boolean variable \({\color{blue} \mathtt{canWriteRHS}_{p}}\).
(See \Cref{alg:AsyncAdditiveRAS}.)
When \({\color{blue} \mathtt{canWriteRHS}_{p}}\) is \texttt{True}, right-hand side data is written to the coarse grid, otherwise this operation is omitted.
Here, the subscripts are used to signify the MPI rank owning the accessed memory region.
As before, index \(0\) corresponds to the (potentially distributed) coarse grid and indices \(1,\dots,P\) correspond to the subdomains.
To improve readability, we show access to a memory region on the calling process in blue, while remote access is printed in red.

In a similar fashion, the coarse grid checks whether every subdomain has written a right-hand side to \({\color{blue} \vec{r}_{0}}\) by polling the state of the local Boolean array \({\color{blue} \mathtt{RHSisReady}_{0}}\).
The communication of the obtained coarse-grid solution back to the subdomains follows the same pattern, using the variables {\color{blue} \(\mathtt{solutionIsReady}_{p}\)}.
The subdomains update their current iterate using the local subdomain solution and the coarse-grid solution.
If the latter is not available, the subdomain solution is used unweighted.
If both solutions are available, then the same weighting \(\left(1/2,1/2\right)\) as in the synchronous case \eqref{eq:additiveRAS} is used.
We note that the algorithm is asynchronous despite the data dependencies.
Coarse grid and subdomain solves do not wait for each other.

We determined by experiments that overall performance is adversely affected if the coarse grid constantly polls the status variable \({\color{blue} \mathtt{RHSisReady}_{0}}\), waiting for all subdomains to provide right-hand side information.
Therefore, we added a sleep statement into its work loop.
If the sleep interval is too short, the sleep statement is ineffective.
If the sleep interval is too large, the coarse grid will be under-used.
Keeping the ratio of attempted coarse-grid solves (i.e. reads from \({\color{blue} \mathtt{RHSisReady}_{0}}\)) to actual performed coarse-grid solves at around \(1/20\) has been proven effective to us.
This can easily be achieved by an adaptive procedure that counts both successful solves and solve attempts and then either increases or decreases the sleep interval accordingly.

\begin{algorithm}[!t]
  \begin{algorithmic}[1]
    \newcommand{\Subdomain}{\STATE \textbf{On subdomains} \begin{ALC@g}}
    \newcommand{\EndSubdomain}{\end{ALC@g}}
    \newcommand{\Coarse}{\STATE \textbf{On coarse grid} \begin{ALC@g}}
    \newcommand{\EndCoarse}{\end{ALC@g}}
    \newcommand{\IND}{\hspace{\algorithmicindent}}
    \newcommand{\INDSTATE}[1][1]{\STATE\hspace{#1\algorithmicindent}}
    \WHILE {not converged}
    \Subdomain
    {
      \STATE Local residual: \(\vec{s}_{p}\leftarrow\mat{D}_{p}\mat{R}_{p}\vec{f} - \mat{A}_{p}\mat{D}_{p}\vec{w}_{p}\)
      \IF{\({\color{blue} \mathtt{canWriteRHS}_{p}}\)}
      {
        \STATE \({\color{red} \vec{r}_{0}} \leftarrow {\color{red}\vec{r}_{0}} + \mat{R}_{0}\mat{R}_{p}^{T}\vec{s}_{p}\)
        \STATE \({\color{blue} \mathtt{canWriteRHS}_{p}} \leftarrow \False\)
        \STATE \({\color{red} \mathtt{RHSisReady}_{0}[p]} \leftarrow \True\)
      }
      \ENDIF
      \STATE Accumulate asynchronously: \(\vec{r}_{p}\leftarrow\sum_{q=1}^{P}\mat{R}_{p}\mat{R}_{q}^{T}\vec{s}_{q}\)
      \STATE Solve: \(\mat{A}_{p}\vec{v}_{p}=\vec{r}_{p}\)
      \IF{{\color{blue} \(\mathtt{solutionIsReady}_{p}\)}}
      {
        \STATE Update: \(\vec{w}_{p}\leftarrow \vec{w}_{p} + \frac{1}{2}\vec{v}_{p} + \frac{1}{2}{\color{blue} \vec{c}_{p}}\)
        \STATE \({\color{blue} \mathtt{solutionIsReady}_{p}} \leftarrow \False\)
      }
      \ELSE
      {
        \STATE Update: \(\vec{w}_{p}\leftarrow \vec{w}_{p}+ \vec{v}_{p}\)
      }
      \ENDIF
    }
    \EndSubdomain
    \Coarse
    {
      \IF{{\color{blue} \(\mathtt{RHSisReady}_{0}[p]\)} \(\forall p=1,\dots,P\)}
      {
        \STATE Solve \(\mat{A}_{0}\vec{v}_{0}={\color{blue} \vec{r}_{0}}\)
        \FOR {\(p=1,\dots,P\)}
        {
          \STATE \({\color{blue} \mathtt{RHSisReady}_{0}[p]} \leftarrow \False\)
          \STATE \({\color{red} \mathtt{canWriteRHS}_{p}} \leftarrow \True\)
          \STATE \({\color{red} \vec{c}_{p}} \leftarrow \mat{R}_{p}\mat{R}_{0}^{T}\vec{v}_{0}\)
          \STATE \({\color{red} \mathtt{solutionIsReady}_{p}} \leftarrow \True\)
        }
        \ENDFOR
      }
      \ELSE
      \STATE Sleep (time adjusted adaptively)
      \ENDIF
    }
    \EndCoarse
    \ENDWHILE
    \Subdomain
    \STATE Post-process synchronously \(\vec{u}_{p}\leftarrow\sum_{q=1}^{P}\mat{R}_{p}\mat{R}_{q}^{T}\mat{D}_{q}\vec{w}_{q}\)
    \EndSubdomain
  \end{algorithmic}
  \caption{
    Asynchronous RAS with additive coarse grid in local form.
    Variables printed in {\color{blue}blue} are exposed memory regions that are local to the calling process.
    {\color{red}Red} variables are remote memory regions.
    Subscripts denote the owning process of the variable.
    Array access is denoted by ``\([\cdot]\)''.
  }\label{alg:AsyncAdditiveRAS}
\end{algorithm}


\subsection{Convergence Analysis of Asynchronous Iterations}
\label{sec:conv-analys-asynchr}


We present below the mathematical framework used to describe and study asynchronous algorithms.
We modify the model introduced by Bertsekas~\cite{Bertsekas1983_DistributedAsynchronousComputationFixedPoints},~\cite{BertsekasTsitsiklis1989_ParallelDistributedComputation} to take into account the fact that data available at a process \(p\) from another process \(q\) might have been produced during different local iterations.
This issue can arise when data is accessed on process \(p\) while it is being overwritten by a new transmission from process \(q\). 

For a mathematical model of these asynchronous iterations on $P$ processors, let us denote by $\left\{ \steering_{n} \right\}_{n\in\N}$ the sequence of non-empty subsets of $\left\{ 1,\dots,P \right\}$, defining which processes update their components at the ``iteration'' $n$, where here ``iteration'' can be thought of as a time stamp.
We call these \emph{sets of update indices}.
Define further for $p,q\in\left\{ 1,\dots,P \right\}$, $\left\{ \vec{\tau}_{q,n}^{(p)}\right\}_{n\in\N}$ a sequence of integer vectors, where $\left(\tau_{q,n}^{(p)}\right)_{i}$, \(1\leq i\leq N_{q}\) represents the iteration number (or time stamp) of the \(i\)-th component of data coming from process $q$ and available on process $p$ at the beginning of the computation of the process which produces $\vec{u}_{p,n}$ at time $n$.
Thus, these are the time stamps of previous computations that are used by process $p$, and thus, the quantities $n - \tau_{q,n,i}^{(p)}$ are sometimes called \emph{delays}.
We use the notation
\begin{align*}
  X_{p}&= \mathbb{R}^{N_{p}}, &\text{and}&   &\Xt &= X_{1}\times\cdots\times X_{P}
\end{align*}
to denote local and global solution spaces, and \(\T_{p,n}:\Xt \rightarrow X_{p}\) the rule that is used to update the local iterate \(\vec{u}_{p,n}\) at iteration \(n\).
We can now define, for each process $p$, the asynchronous iterations as follows:
\begin{align}
  \vec{u}_{p,n} &=
  \begin{cases}
    \T_{p,n}\left(
      \vec{u}_{1,n}^{(p)},\dots,
      \vec{u}_{P,n}^{(p)}\right) &\text{if } p\in \steering_{n}, \\
    \vec{u}_{p,n-1}&\text{if } p \notin \steering_{n}.
  \end{cases}
  \label{eqn:AsynchronousIterativeAlgorithm}
\end{align}
The iteration is initialized using some initial guess for \(\vec{u}_{p,0}\), and we used the notation \(\vec{u}_{q,n}^{(p)}:=\vec{u}_{q,\tau_{q}^{(p)}(n)}\) to denote the data from process \(q\) that is available to process \(p\) at time \(n\).

In other words, at time $n$, either $u_{p,\bullet}$ is not updated (if $p \notin \steering_{n}$) or it is updated with the result of applying the (local) operator $\T_{p,n}$ to the variables computed at times $\vec{\tau}_{\bullet}^{(p)}$.
For comparison, the corresponding synchronous iteration is given by
\begin{align}
  \vec{u}_{p,n} &= \T_{p,n}\left( \vec{u}_{1,n-1},\dots, \vec{u}_{P,n-1}\right),
  \label{eqn:SynchronousIterativeAlgorithm}
\end{align}
or, in compact form, as
\begin{align}
      \tilde{u}_{n}&=\Tt_{n}\left(\tilde{u}_{n-1}\right),
\end{align}
where
\begin{align*}
  \tilde{u}_{n} &= \left( \vec{u}_{1,n},\dots, \vec{u}_{P,n}  \right), &\text{and}&& \Tt_{n} &= \left( \T_{1,n}, \T_{2,n},\cdots,\T_{P,n} \right).
\end{align*}

We further assume that the three following conditions are satisfied
\begin{eqnarray}
  \forall p,q\in\left\{ 1,\dots,P\right\}, 1\leq i\leq N_{q}, \forall n\in\N^*, \left(\tau_{q,n}^{(p)}\right)_{i} \leq n ,
  \label{eqn:PositiveRetarded}
\\
  \forall p\in\left\{ 1,\dots,P \right\}, \operatorname{card}\left\{ n\in\N^* \mid p\in \steering_{n} \right\} = \infty ,
  \label{eqn:NoStop}
\\\
  \forall p,q\in\left\{ 1,\dots,P \right\}, 1\leq i\leq N_{q}, \lim_{n\rightarrow+\infty} \left(\tau_{q,n}^{(p)}\right)_{i} = \infty .
 \label{eqn:FiniteDelay}
\end{eqnarray}
Condition~\eqref{eqn:PositiveRetarded} indicates that data used at the time $n$ must have been produced before time $n$, i.e., time does not flow backward.
Condition~\eqref{eqn:NoStop} means that no process will ever stop updating its components.
Condition~\eqref{eqn:FiniteDelay} corresponds to the fact that new data will always be provided to the process.
In other words, no process will have a piece of data that is never updated.

We note that these assumptions pose no significant restrictions on the iterations that we consider, but are necessary for the analysis.

Assume that each \(X_{p}\) is a normed linear space, equipped with a norm \(\norm{\cdot}_{p}\).
Given a positive vector \(w\in\mathbb{R}_{>0}^{P}\), the weighted norm \(\norm{\cdot}_{w}\) on the product space \(X\) is defined to be
\begin{align*}
  \norm{\tilde{u}}_{w}&=\max_{p=1,\dots,P}\frac{\norm{\vec{u}_{p}}_{p}}{w_{p}}\cdot
\end{align*}

We are ready to present a convergence theorem for asynchronous iterative algorithms, whose proof can be found in~\cite[Theorem 3.3]{FrommerSzyld2000_AsynchronousIterations}.
\begin{theorem}
  Assume that there exists \(\tilde{u}^{*}\in X\) such that \(\Tt_{n}\left(\tilde{u}^{*}\right)=\tilde{u}^{*}\) for all \(n\).
  Moreover, assume that there exists \(\gamma\in[0,1)\) and \(w\in\mathbb{R}^{P}_{>0}\) such that for all \(n\) we have
  \begin{align*}
    \norm{\Tt_{n}\left(\tilde{u}\right)-\tilde{u}^{*}}_{w} &\leq\gamma \norm{\tilde{u}-\tilde{u}^{*}}_{w}.
  \end{align*}
  Then the asynchronous iterates \(\tilde{u}_{n}\) converge to \(\tilde{u}^{*}\), the unique common fixed point of all \(\Tt_{n}\).
\end{theorem}

In view of equations~\eqref{eq:4} and~\eqref{eq:5}, we have
\begin{align*}
  \T_{p,n}^{1L}(\vec{w}_{1},\dots,\vec{w}_{P})&=\vec{w}_{p} + \mat{A}_{p}^{-1}\sum_{q=1}^{P} \mat{R}_{p}\mat{R}_{q}^{T}\left(\mat{D}_{q}\mat{R}_{q}\vec{f} - \mat{A}_{q}\mat{D}_{q}\vec{w}_{q}\right)
\end{align*}
for the one-level method.
We immediately observe that the mappings \(\T_{p,\bullet}^{1L}\) do not depend on \(n\), and that the iteration is stationary.

In order to tackle the two-level method, based on \Cref{alg:AsyncAdditiveRAS} we set
\begin{align*}
  \T_{0,n}^{2L}(\vec{v}_{0},\vec{w}_{1},\dots,\vec{w}_{P})&=\mat{A}_{0}^{-1}\sum_{q=1}^{P} \mat{R}_{0}\mat{R}_{q}^{T}\left(\mat{D}_{q}\mat{R}_{q}\vec{f} - \mat{A}_{q}\mat{D}_{q}\vec{w}_{q}\right).
\end{align*}
Moreover, for \(p=1,\dots,P\), we set 
\begin{align*}
  \T_{p,n}^{2L}(\vec{v}_{0},\vec{w}_{1},\dots,\vec{w}_{P})&=\vec{w}_{p} + \frac{1}{2}\mat{R}_{p}\mat{R}_{0}^{T}\vec{v}_{0} + \frac{1}{2}\mat{A}_{p}^{-1}\sum_{q=1}^{P} \mat{R}_{p}\mat{R}_{q}^{T}\left(\mat{D}_{q}\mat{R}_{q}\vec{f} - \mat{A}_{q}\mat{D}_{q}\vec{w}_{q}\right)
\end{align*}
for iteration numbers \(n\) that include coarse-grid updates, and
\begin{align*}
  \T_{p,n}^{2L}(\vec{v}_{0},\vec{w}_{1},\dots,\vec{w}_{P})&=\T_{p,n}^{1L}(\vec{w}_{1},\dots,\vec{w}_{P})
\end{align*}
for iterations \(n\) without coarse-grid update.
Status variables such as \(\mathtt{canWriteRHS}_{p}\) act implicitly as constraints on the sets of update indices \(\steering_{n}\) and do not appear in the definition of the mappings \(\T_{p,n}^{2L}\).

It has been shown in \cite{FrommerSzyld2001_AlgebraicConvergenceTheoryRestricted} that both the one- and the two-level iterations are contracting in a weighted max-norm, provided that \(\mat{A}\) is a non-singular M-matrix, i.e. if \(\mat{A}\) has nonpositive off-diagonal elements and all entries of \(\mat{A}^{-1}\) are nonnegative.

Thus, we have the following result.
\begin{theorem}\label{thm:async-conv}
  The one-level method given in~\eqref{eq:4} and~\eqref{eq:5} and the two-level method given in~\Cref{alg:AsyncAdditiveRAS} converge, provided that \(\mat{A}\) is a non-singular M-matrix and that the conditions~\ref{eqn:PositiveRetarded}--\ref{eqn:FiniteDelay} hold. 
\end{theorem}

For further extensions of the theory, such as inexact sub-solves with \(\mat{A}_{p}^{-1}\) replaced by some (potentially nonstationary) \(\mat{S}_{p,n}\approx \mat{A}_{p}^{-1}\), we refer the reader to \cite{FrommerSzyld2001_AlgebraicConvergenceTheoryRestricted}. 


\section{One-sided Message Passing Interface} \label{sec:one-sided-mpi}
In order to drive the asynchronous method in a distributed memory setting, we use a one-sided approach wherein the remote process incurs minimal overhead for servicing received messages from the sender process.
The one-sided approach is achieved in MPI using the Remote Memory Access (RMA) semantics, wherein every process exposes a part of its local memory window to remote processes for read as well as write operations.
However, in reality, a synchronization between the source and the target process is required for progress of the underlying application.
This active synchronization step, while still preserving the asynchronous nature of the algorithm, is expensive and might erode the natural gains  obtained from the asynchronous method.
Therefore in order to extract the maximum gains from an asynchronous method, a passive approach is required.
A passive approach entails transmission of messages which causes little to no interference to the target process.
As a result, the target process does not need to yield its operating system time for servicing incoming message interrupts and therefore does not participate in the communication process.
The RMA framework on MPI implements passive target synchronization with the help of two sets of primitives \mpi{Win_lock}/\mpi{Win_unlock} and \mpi{Win_lock_all}/\mpi{Win_unlock_all}.
While the former involves opening and closing the exposure epoch on remote nodes for each access operation, the latter only requires opening and closing of access epoch once during the application lifetime incurring less target synchronization overhead.


RMA's passive one-sided communication can leverage a hardware mechanism known as Remote Direct Memory Access (RDMA) \cite{LiuJiuxingEtAl2004_RDMA} when available.
It allows RMA to directly map memory windows to the RDMA engine, allowing messages written by remote processes to be directly read by each process.
This leads to minimum disturbance to the remote process and achieving a truly passive, one-sided communication scheme.


RDMA is usually a hardware characteristic that may not be supported by all machines.
Though we expect one-sided communication of RMA to be able to handle progress of communication in an entirely asynchronous manner, it generally fails to do so since MPI does not guarantee asynchronous progress.
In such a case, asynchronous progress may be enforced by allocating certain auxiliary cores to ghost processes that solely perform the task of asynchronous progress control.
As a consequence we obtain an RDMA agnostic system while simultaneously obtaining the benefits of RDMA.
Even in the presence of RDMA, asynchronous progress control mechanism can be complementary since the low level RDMA engine may not be capable to handle high volumes of communication.
Casper \cite{SiPenaEtAl2015_Casper} and Intel Asynchronous Progress Control (APC) are two such implementations that provide ghost processes for asynchronous progress control.



\section{Implementation and numerical experiments}
\label{sec:experiments}

\subsection{Comparison of one-sided MPI communication options}
\label{sec:selection-one-sided}

There are a multitude of options for achieving asynchronous neighborhood exchange.
Data that is supposed to be moved from rank \(p\) to rank \(q\) could be held in MPI windows on either \(p\) or \(q\).
In the first case, rank \(p\) will write the data to its local buffer using \mpi{Put}, and rank \(q\) will retrieve it from the remote buffer using \mpi{Get}.
In the second case, rank \(p\) writes the data to the remote memory region using a \mpi{Put}, and \(q\) retrieves using a local \mpi{Get}.

The second distinction comes from the type of locking mechanism used.
Exclusive or shared locks can either be applied for each individual memory access (\mpi{Win_lock}), or windows can be locked in shared fashion for all subsequent access (\mpi{WIN_lock_all}).
In the latter case, windows can be flushed using any of the available flush operations.

We benchmark the different available options in a simple test case in order to determine which one should be used in the implementation of our domain decomposition solvers.
The performance of one-sided MPI communication depends on the support provided by the MPI implementation as well as the network hardware.
These experiments are performed on the Haswell partition of Cori at the National Energy Research Scientific Computing Center (NERSC), using the default Cray MPICH, version 7.7.3.
Since one-sided MPI has not been widely adopted, performance variations compared to the classical two-sided routines can be expected to be much more significant.
It should be noted that different network hardware and better support in future MPI versions could further improve timings for one-sided MPI routines.

64 MPI ranks are arranged in a three dimensional regular periodic grid (3D torus), and each rank repeatedly exchanges a vector of doubles with its \(26\) neighbors.
This test mimics the communication pattern in the neighborhood exchange of the one-level method.
For each of the possible communication option as given in \Cref{tab:comm}, we measure the time it takes to perform 50,000 exchanges of vectors of 500 doubles.
By exchanging vectors that have a constant value corresponding to the exchange iteration, we can also measure how often inconsistent data is accessed (i.e. data that is accessed before it has been completely been transmitted).
This phenomenon does not occur when using two-sided communication, since completion is guaranteed by the implementation.
While the absolute number of accesses to incomplete writes is probably quite dependent on the ratio of computation to communication, we are interested in the susceptibility of the different communication options.

\begin{table}
  \centering
  {\renewcommand\arraystretch{1.5}
    \begin{tabular}{p{0.07\linewidth}p{0.35\linewidth}p{0.2\linewidth}p{0.09\linewidth}p{0.1\linewidth}}
      global lock & per comm phase & per neighbor & time in seconds & inconsistency fraction\\
      \hline
      \xmark &  \mpi{Win_lock}(EXCLUSIVE) \newline \mpi{Win_unlock} & local \mpi{Put}, \newline remote \mpi{Get} &  34.6 & 0.0 \\ \hline
      \xmark &  \mpi{Win_lock}(EXCLUSIVE) \newline \mpi{Win_unlock} & remote \mpi{Put}, \newline local \mpi{Get} &  37.8 & 0.0 \\ \hline
      \xmark &  \mpi{Win_lock}(SHARED) \newline \mpi{Win_unlock}  & local \mpi{Put}, \newline remote \mpi{Get} &  31.8 & 0.00151 \\ \hline
      \xmark &  \mpi{Win_lock}(SHARED) \newline \mpi{Win_unlock}  & remote \mpi{Put}, \newline local \mpi{Get} &  33.0 & 0.00254 \\ \hline
      n/a & \mpi{Wait_all} & \mpi{Isend}, \newline \mpi{Irecv} &  9.59 & 0.0 \\ \hline
      \cmark &  - & local \mpi{Put}, \newline remote \mpi{Get} &  25.8 & 0.123 \\ \hline
      \cmark &  - & remote \mpi{Put}, \newline local \mpi{Get} &  8.42 & 0.00716 \\ \hline
      \cmark &  \mpi{flush_all}& local \mpi{Put}, \newline remote \mpi{Get} &  22.1 & 0.117 \\ \hline
      \cmark &  \mpi{flush_all}& remote \mpi{Put}, \newline local \mpi{Get} &  9.06 & 0.00491 \\ \hline
      \cmark &  \mpi{flush_local_all}& local \mpi{Put}, \newline remote \mpi{Get} &  22.1 & 0.099 \\ \hline
      \cmark &  \mpi{flush_local_all}& remote \mpi{Put}, \newline local \mpi{Get} &  9.02 & 0.00501 \\ \hline
      \cmark &  \mpi{flush_local}& local \mpi{Put}, \newline remote \mpi{Get} &  24.1 & 0.172 \\ \hline
      \cmark &  \mpi{flush_local}& remote \mpi{Put}, \newline local \mpi{Get} &  10.7 & 0.00198 \\ \hline
      \cmark &  \mpi{flush} & local \mpi{Put}, \newline remote \mpi{Get} &  21.8 & 0.105 \\ \hline
      \cmark &  \mpi{flush} & remote \mpi{Put}, \newline local \mpi{Get} &  11.2 & 0.00207
    \end{tabular}
  }
  \caption{
    Results of communication test described in \Cref{sec:selection-one-sided} on 64 MPI ranks.
    The listed operations are either performed once per neighborhood communication phase, or for each individual neighborhood exchange.
    If \texttt{MPI\_Win\_lock\_all}/\texttt{MPI\_Win\_unlock\_all} is used, the column ``global lock'' has a \cmark.
    We measured the time for 50,000 repetitions and the fraction of neighborhood exchanges leading to incompletely written data.
  }
  \label{tab:comm}
\end{table}

We make several observations.
Unsurprisingly, the use of exclusive locks does not perform well in terms of time.
However, the use of shared locks in every communication phase performs equally poorly, which is why we decide to use global locking and unlocking (\mpi{Win_lock_all} / \mpi{Win_unlock_all}) in what follows.
Using global locking, we see that using remote puts instead of remote gets is significantly faster.

We also observe that unless exclusive locks are used, we always experience access to inconsistent data.
This might not be of too much importance within our application, since it amounts to using residual information that is only slightly more outdated.
Finally, we observe that using global locking and puts results in faster communication than classical two-sided non-blocking communication.

Based on the above results, we choose to use global locking using \mpi{Win_lock_all} / \mpi{Win_unlock_all}, paired with remote \mpi{Put}s and local \mpi{Get}s and \mpi{flush_all}, since it appears to provide a good balance of speed and consistency.
We note however that these results might depend significantly on characteristics of the system and the MPI implementation.


\subsection{Performance metrics}

The average contraction factor per iteration is defined as \(\tilde{\rho} = \left({r_{\text{final}}}/{r_{0}}\right)^{\frac{1}{K}}\), where \(r_{0}\) is the norm of the initial residual vector, \(r_{\text{final}}\) the norm of the final residual vector, and \(K\) is the number of iterations that were taken to decrease the residual from \(r_{0}\) to \(r_{\text{final}}\).
For an asynchronous method, the number of iterations varies from subdomain to subdomain, and hence \(\tilde{\rho}\) is not well-defined.
The following generalization permits us to compare synchronous methods with their asynchronous counterpart:
\begin{align*}
  \hat{\rho} = \left(\frac{r_{\text{final}}}{r_{0}}\right)^{\frac{\tau_{\text{sync}}}{T}}.
\end{align*}
Here, \(T\) is the total iteration time, and \(\tau_{\text{sync}}\) is the average time for a single iteration of the synchronous method.
In the synchronous case, since \(T=\tau_{\text{sync}} K\), \(\hat{\rho}\) recovers \(\tilde{\rho}\).
The approximate contraction factor \(\hat{\rho}\) can be interpreted as the average contraction of the residual norm in the time of a single synchronous iteration.
As it will be visible in the results to follow, we note here that \(\hat{\rho}\) for the asynchronous method obviously depends on the total iteration time for the synchronous method.
Assume that the total iteration time for the synchronous method doubles, but the time taken by the asynchronous one stays constant.
Consequently, the approximate contraction factor for the synchronous method stays constant, but the contraction factor for the asynchronous method gets squared and therefore decreases.

\subsection{Test problem}
\label{sec:test-problem}

As a test problem, we solve
\begin{align*}
  -\Delta u&=f & \text{in }\Omega = [0,1]^{d}, && u&=0 & \text{on }\partial\Omega,
\end{align*}
where the right-hand side is \(f=d\pi^{2}\prod_{k=1}^{d}\sin\left(\pi x_{k}\right)\).
The corresponding solutions is \(u=\prod_{k=1}^{d}\sin\left(\pi x\right)\).
We discretize \(\Omega\) using a uniform simplicial mesh and approximate the solution using piece-wise linear finite elements.
We note that the arising system matrix \(\mat{A}\) is a non-singular M-matrix, and therefore \Cref{thm:async-conv} applies.
Furthermore, we mention that the generalization of the test problem to convection-diffusion problems with non-constant diffusion coefficient is possible, but does not alter the numerical results obtained below in a significant way, which is why we only present the case of the standard Poisson problem.

\subsection{Convergence detection}
\label{sec:conv-detection}

In classical synchronous iterative methods, a stopping criterion of the form \(\norm{\vec{r}}<\varepsilon\) is evaluated at every iteration.
Here, \(\vec{r}\) is the residual vector, \(\varepsilon\) is a prescribed tolerance (that might be chosen as a function of the discretization error), and \(\norm{\cdot}\) is an appropriate norm.
The global quantity \(\norm{r}\) needs to be computed as the sum of local contributions from all the subdomains.
This implies that convergence detection in asynchronous methods is not straightforward, since collective communication primitives require synchronization.
In the numerical examples below, we use a simplistic convergence criterion, consisting in writing the local contributions to a master rank, say rank 0.
This master rank sums the contributions, and determines if this approximation of the global residual norm is smaller than the prescribed tolerance.
If so, the master rank declares global convergence and notifies the other ranks by sending a non-blocking message. 
This simplistic convergence detection mechanism has several drawbacks.
For one, the global residual is updated by the master rank, which might not happen frequently enough.
Hence it is possible that the iteration continues despite the true global residual norm already being smaller than the tolerance.
Moreover, the mechanism puts an increased load on the network connection to the master rank, since every subdomain writes to its memory region.
Finally, since the local contributions to the residual norm are not necessarily monotonically decreasing, the criterion might actually detect convergence when the true global residual is not yet smaller than the tolerance.
The delicate topic of asynchronous convergence detection has been treated in much detail in the literature, and we refer to \cite{BahiContassot-VivierEtAl2005_DecentralizedConvergenceDetectionAlgorithm,MagoulesGbikpi-Benissan2017_DistributedConvergenceDetectionBased} for an overview of more elaborate approaches.
While these detection schemes mostly address the shortcomings of the above approach, their correct implementation turns out to be quite involved.
Since we are not observing any major issues with our simplistic convergence detection scheme for the test problems that we consider, we have not implemented any of the schemes available in the literature.

\subsection{Platform and implementation details}
\label{sec:platf-impl-deta}

All runs are performed on the Haswell partition of the Cori supercomputer at NERSC.
While all of the code was written from scratch, the differences between the synchronous and the asynchronous code path are limited, since only the communication layer and the stopping criterion need to be changed.
(E.g. compare \Cref{alg:additiveRAS,alg:AsyncAdditiveRAS}.)
We stress that the asynchronous solver uses one-sided communication only in the solve phase.
Therefore, we record solve times only, since the time to set up the solver is unaffected by the type of communication in the solution phase.
Furthermore, all subdomains are synchronized via a \mpi{Barrier} before entering the solve phase.
One MPI rank is used per core, i.e. 32 ranks per Haswell node.
Moreover, one subdomain is assigned to each MPI rank.
The underlying mesh is partitioned either into uniformly sized rectangular subdomains or using the METIS library \cite{KarypisKumar1998_FastHighQualityMultilevel}.
In the latter case, the option to minimize the overall communication volume is used.
Our solvers handle general unstructured matrices, and the structure of the mesh is not exploited.
We either use
\begin{itemize}
\item direct solvers for subdomain and coarse-grid problems, provided by SuperLU \cite{LiDemmelEtAl1999_SuperluUsersGuide,DemmelEisenstatEtAl1999_SupernodalApproachToSparsePartialPivoting}, or
\item conjugate gradient method preconditioned with an incomplete Cholesky factorization for the subdomain problems and a geometric multigrid solver for the coarse problem.
\end{itemize}
The latter option would allow for a distributed coarse-grid solve and is therefore in principle more scalable.
In all numerical examples, we will use only a single core for the coarse-grid solve.

\subsection{Comparison against HPDDM}
\label{sec:comp-against-hpddm}

\begin{figure}[!t]
  \centering
  \includegraphics{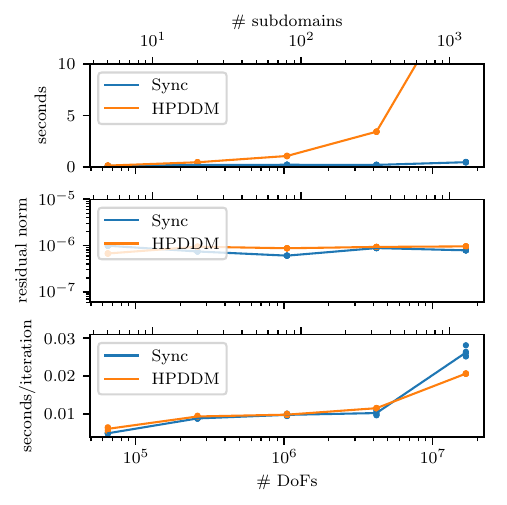}
  \caption{
    Weak scaling of GMRES preconditioned by two-level additive RAS using the synchronous version of our code (Sync) and HPDDM for the 2D test problem in a load balanced case.
    From top to bottom: Total solution time, final residual norm, and time per iteration.
    Mean values are given by solid lines, and individual runs as dots.
  }
  \label{fig:weakScaling2DHPDDM}
\end{figure}

\begin{figure}[!t]
  \centering
  \includegraphics{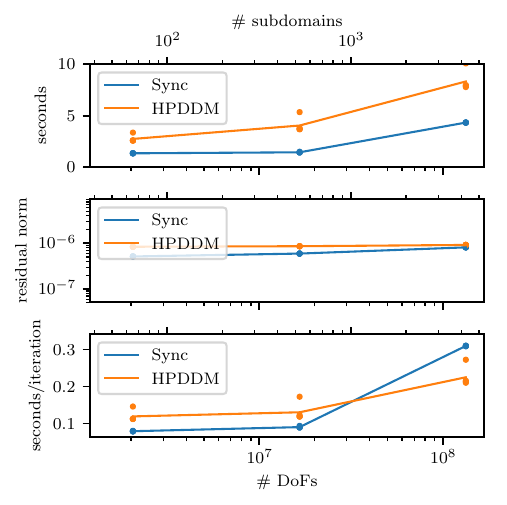}
  \caption{
    Weak scaling of GMRES preconditioned by two-level additive RAS using the synchronous version of our code (Sync) and HPDDM for the 3D test problem in a load balanced case.
    From top to bottom: Total solution time, final residual norm, and time per iteration.
    Mean values are given by solid lines, and individual runs as dots.
  }
  \label{fig:weakScaling3DHPDDM}
\end{figure}

We verify the performance of the synchronous version of our code against the HPDDM library \cite{hpddm,JolivetHechtEtAl2013_ScalableDomainDecompositionPreconditioners,DoleanJolivetEtAl2015_IntroductionToDomainDecompositionMethods} using the 2D and 3D test problems from \Cref{sec:test-problem}.
In all cases, we set up a GMRES solver and use a two-level additive RAS as right preconditioner.
The reason for using a Krylov method here is that domain decomposition methods in general are commonly used as preconditioners, and HPDDM is most likely developed with that use case in mind.
HPDDM was linked against the Intel Math Kernel Library, SuiteSparse~\cite{ChenDavisEtAl2008_Algorithm887} and ARPACK~\cite{LehoucqSorensenEtAl1998_ArpackUsersGuide} and the option for coarse grid data exchange using \mpi{Igather}/\mpi{Iscatter} was enabled.
We use the following parameters in HPDDM: \texttt{-hpddm\_krylov\_method=gmres -hpddm\_schwarz\_method=ras -hpddm\_schwarz\_coarse\_correction=additive \newline -hpddm\_geneo\_nu=NU}, where \texttt{NU} is chosen so that the size of the coarse grid matches our solver.
In 2D the subdomains consist of roughly 20k unknowns, and the coarse grid contains about 16 unknowns per subdomain.
In 3D the subdomains consist of roughly 40k unknowns, and the coarse grid contains about 1 unknown per subdomain.
In \Cref{fig:weakScaling2DHPDDM,fig:weakScaling3DHPDDM} we plot the results of weak scaling experiments: overall solve time, the reached residual norm and the time per iteration.
We repeated each run 5 times.
Mean values are given by solid lines, and individual runs as dots.
We observe that while the time to convergence behaves quite differently for both implementations, the time per iteration follows the same trend.
HPDDM behaves slightly better at large subdomain count which could be explained by the use of MPI collectives, but might also be an artifact of the difference in convergence behavior or the difference in coarse solvers (HPDDM uses Cholmod).
We can therefore use our synchronous method as a base of comparison for the newly developed asynchronous solver.

\subsection{One-level RAS, 2D test problem, strong scaling}
\label{sec:one-level-ras}

\begin{figure}[!t]
  \centering
  \includegraphics{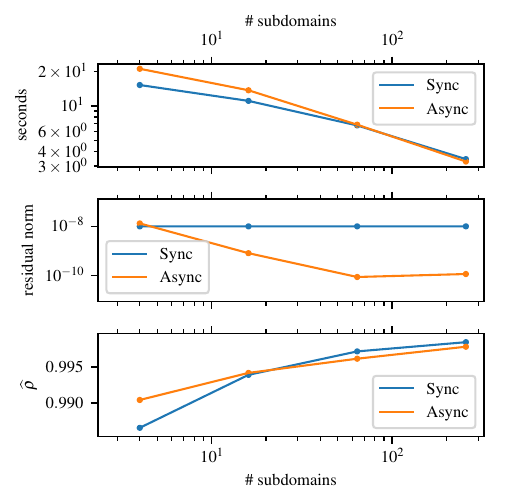}
  \caption{
    Strong scaling of synchronous and asynchronous one-level RAS for the 2D test problem with system size of approximately 261,000 unknowns.
    The subdomains are load balanced.
    From top to bottom: Solution time, final residual norm, and the resulting approximate contraction factor \(\hat{\rho}\).
    It can be observed that the synchronous method is significantly faster than for smaller numbers of subdomains (cores), yet comparatively slower for larger number of subdomains, as shown by the contraction factor.
  }
  \label{fig:strongScaling1level}
\end{figure}

We compare synchronous and asynchronous one-level RAS in a strong scaling experiment, where we fix the global problem size of a 2D test problem to about 261,000 unknowns, and vary the number of subdomains between 4 and 256.
We cannot expect good scaling behavior for this one level method, since increasing the number of subdomains adversely affects the rate of convergence.
The iteration is terminated based on the simplistic convergence criterion described in \Cref{sec:conv-detection}.
In \Cref{fig:strongScaling1level} we show solve time, final residual norm and approximate rate of convergence.
It can be observed that the synchronous method is faster for smaller numbers of subdomains, yet comparatively slower for larger number of subdomains.
The crossover point between the two regimes appears to be at 64 subdomains.

\begin{figure}[!t]
  \centering
  \includegraphics{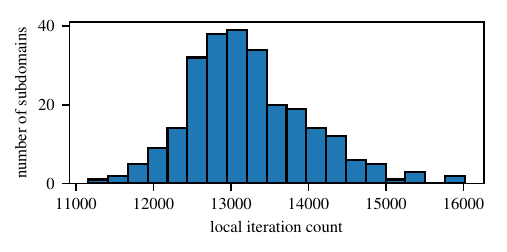}
  \caption{
    Histogram of local iteration counts for asynchronous one-level RAS for the 2D test problem with 256 subdomains (load balanced case).
  }
  \label{fig:iterationHistogram}
\end{figure}

An important question is whether the asynchronous method happens to converge because every subdomain performs the same number of local iterations, and hence the asynchronous method just mirrors the synchronous one, merely with a different communication method.
The histogram in \Cref{fig:iterationHistogram} shows that this is not the case.
The number of local iterations varies significantly.
The slowest subdomain performs barely more than 11,000 iterations, whereas the fastest one almost reaches 16,000.
The problem was load balanced by the number of degrees of freedom in each subdomain, thus the local solves are also approximately balanced but the communication is likely slightly imbalanced.
This means that in this scenario, system and network noise are the main contributions to the observed variations in local iteration counts.
For comparison, when only 4 subdomains were used, the local iteration counts were 1497, 1500, 1504 and 1527.

\setlength{\abovecaptionskip}{-10pt}

\begin{figure}[!t]
  \centering
  \includegraphics{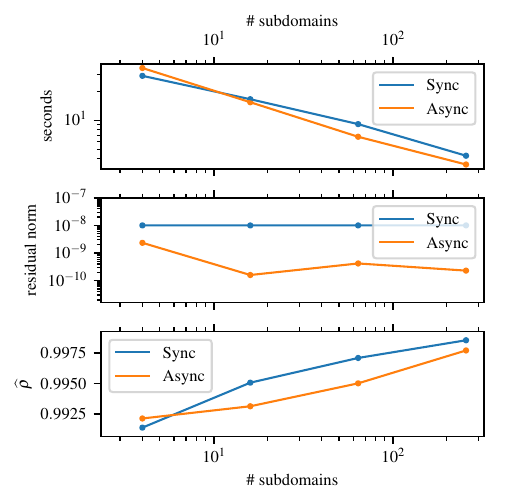}
  \caption{
    Strong scaling of synchronous and asynchronous one-level RAS for the 2D test problem with a system size of approximately 261,000 unknowns under load imbalance: one subdomain is 50\% larger than the rest.
    From top to bottom: Solution time, final residual norm, and approximate contraction factor \(\hat{\rho}\).
    It can be observed that the asynchronous method outperforms the synchronous one in all but the 4 subdomain case, as shown by the contraction factor.
    The advantage of the asynchronous method over the synchronous one is increased, as compared to \Cref{fig:strongScaling1level}.
  }
  \label{fig:strongScaling1levelImbalance}
\end{figure}

The advantage of asynchronous RAS becomes even clearer when the experiment is repeated under load imbalance.
We create an artificial load imbalance by choosing one of the subdomains to be 50\% larger than the rest.
In \Cref{fig:strongScaling1levelImbalance} it is observed that the asynchronous method outperforms the synchronous one in all but the 4 subdomain case.

\subsection{Two-level RAS, 2D test problem}
\label{sec:two-level-ras}

In order to gauge the performance and scalability of the synchronous and asynchronous two-level RAS solvers, we perform weak and strong scaling experiments.

\subsubsection{Weak scaling}

In the weak scaling experiment the number of subdomains \(P\) and the global number of degrees of freedom (DoFs) are increased proportionally.
We use 16, 64, 256 and 1024 subdomains to solve the 2D test problem.
The local number of unknowns on each subdomain is kept constant at almost 20,000.
The coarse-grid problem increases in size proportionally to the number of subdomains, with approximately 16 unknowns per subdomain.
Again, the iteration is terminated based on the simplistic convergence criterion described in \Cref{sec:conv-detection}.

\begin{figure}[!t]
  \centering
  \includegraphics{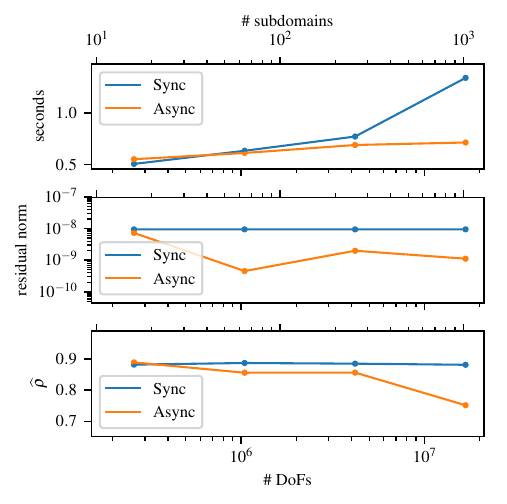}
  \caption{
    Weak scaling of synchronous and asynchronous two-level additive RAS for the 2D test problem, load balanced case.
    From top to bottom: Total solution time, final residual norm, and approximate contraction factor \(\hat{\rho}\).
    One can observe that for 16, 64 and 256 subdomains, the asynchronous and the synchronous method take almost the same time for the solve, with a slight advantage for the asynchronous method.
    For 1024 subdomains, however, the synchronous method is seen to take significantly more time, since the coarse grid, due to its size, starts to be the limiting factor.
    The asynchronous method is not affected by this.
  }
  \label{fig:weakScaling2level}
\end{figure}
\begin{figure}[!t]
  \centering
  \includegraphics{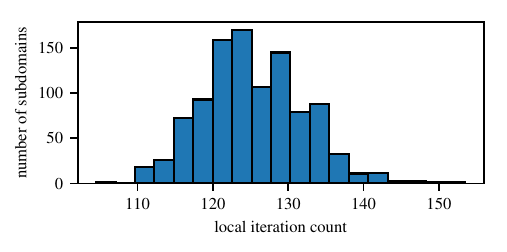}
  \caption{
    Histogram of local iteration counts asynchronous two-level additive RAS for the 2D test problem with 1024 subdomains.
  }
  \label{fig:iterationHistogram2level}
\end{figure}

In \Cref{fig:weakScaling2level} we plot the solution time, the achieved residual norm and the average contraction factor \(\hat{\rho}\) depending on the global problem size.
Both the synchronous and the asynchronous method reach the prescribed tolerance of~\(10^{-8}\).
Due to the lack of an efficient mechanism of convergence detection, the asynchronous method ends up iterating longer than necessary, so that the final residual norm often is smaller than~\(10^{-9}\).
The number of iterations in the synchronous case is about 110, whereas the number of local iterations in the asynchronous case varies between 110 and 150.
(See \Cref{fig:iterationHistogram2level}.)
The iteration counts are significantly lower than for the one-level methods.
One can observe that for 16, 64 and 256 subdomains, the asynchronous and the synchronous methods take almost the same time for the solve.
For 1024 subdomains, however, the synchronous method is seen to take significantly more time.
This can be explained by the fact that for 1024 subdomains, the size of the coarse grid is comparable to the size of the subdomains, and hence the coarse-grid solve which exchanges information with all the subdomains slows down the overall progress.
For the asynchronous case this is not observed, since the subdomains do not have to wait for information from the coarse grid.
This explains why we see better weak scalability for the asynchronous method than for the synchronous variant, and why we can observe a speedup of 2x of the asynchronous method over its synchronous counterpart.
The third subplot of \Cref{fig:weakScaling2level} shows that the asynchronous method outperforms its synchronous equivalent in all but the smallest problem.

\begin{figure}[!t]
  \centering
  \includegraphics{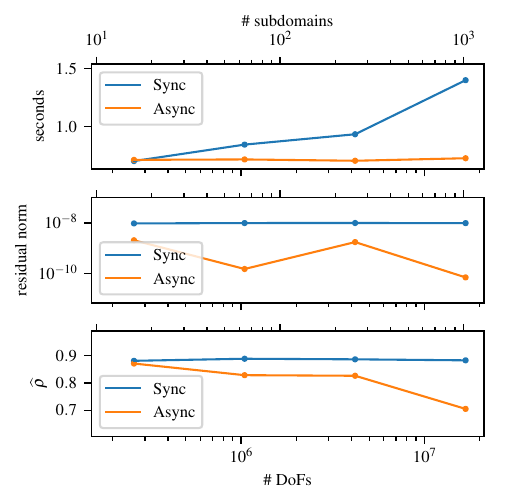}
  \caption{
    Weak scaling of synchronous and asynchronous two-level additive RAS for the 2D test problem under load imbalance: one subdomain is 50\% larger than all the other ones.
    From top to bottom: Total solution time, final residual norm, and approximate contraction factor \(\hat{\rho}\).
    The advantage of the asynchronous method over the synchronous one is increased, as compared to \Cref{fig:weakScaling2level}.
  }
  \label{fig:weakScaling2levelImbalance}
\end{figure}

To further illustrate the effect of load imbalance, we repeat the previous experiment with one subdomain being 50\% larger than the rest.
The results are shown in \Cref{fig:weakScaling2levelImbalance}.
While the results are mostly consistent with the previous case, it can be seen that, as expected, the performance advantage of the asynchronous method over the synchronous one has increased.
Even before the size of the coarse-grid system is comparable to the size of the typical subdomain problem, the asynchronous method outperforms its synchronous counterpart.

\subsubsection{Strong scaling}

For the strong scaling experiment the global number of degrees of freedom used to discretize the 2D test problem is fixed at about 4 million.
The coarse-grid problem consists of approximately 4,000 unknowns.
The number of subdomains used on the fine level takes values in \(\left\{4,16,64,256\right\}\).
This means that the coarse-grid problem is always smaller than the typical subdomain problem, and no slowdown due to an imbalance of the computational cost of coarse and fine solve should arise.

\begin{figure}[!t]
  \centering
  \includegraphics{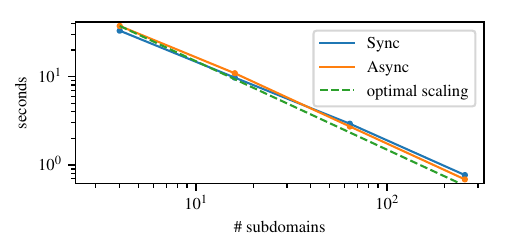}
  \includegraphics{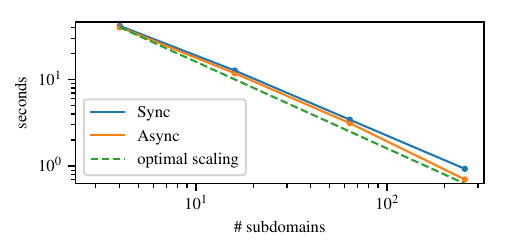}
  \caption{
    Strong scaling of synchronous and asynchronous two-level additive RAS for the 2D test problem.
    On top: load balanced subdomains.
    At the bottom: load imbalance, one subdomain is 50\% larger than the others.
  }
  \label{fig:strongScaling2level}
\end{figure}

The timing results are shown in the top of \Cref{fig:strongScaling2level}.
Both synchronous and asynchronous method display good strong scaling behavior.
It is observed that the synchronous method is faster than the asynchronous method for smaller subdomain count.
But already for 64 subdomains this behavior is reversed, and the asynchronous method outperforms the synchronous one.
This suggests that synchronization is an important factor already at modest core count.

At the bottom of \Cref{fig:strongScaling2level}, we show the timing results in the case of load imbalance.
It can be seen that the asynchronous method is faster than the synchronous one independent of the number of subdomains, and that its performance advantage increases as more processes are used.

\subsection{Two-level RAS with iterative sub-solves, 3D test problem}

The density of the subdomain matrices \(\mat{A}_{p}\) in 3D (about 15 entries per row) is higher than for the 2D test problem (about 7 entries per row).
This means that direct factorization leads to more fill-in and thereby is more expensive.
Therefore, we solve subdomain and coarse problem of the three dimensional test case using iterative solvers.
For the subdomains, we use a conjugate gradient solver preconditioned by an incomplete Cholesky factorization.
We employ a relative tolerance of \(1/10\) which has been determined experimentally to be sufficient.
The coarse-grid problem is solved using a single V-cycle of a geometric multigrid solver with one step of Gauss-Seidel for pre- and post-smoothing.
The use of multigrid allows us to solve the coarse-grid problem in a distributed fashion when it becomes too large for a single MPI rank.
The global problem is partitioned into uniformly sized regular subdomains.
The local number of unknowns on each subdomain is kept constant at about 40,000.
The coarse-grid problem increases in size proportionally to the number of subdomains, with approximately one unknown per subdomain.

\begin{figure}
  \centering
  \includegraphics{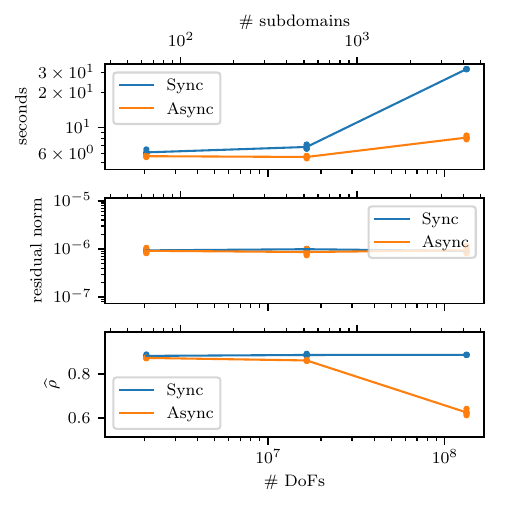}
  \caption{
    Weak scaling of synchronous and asynchronous two-level additive RAS, load balanced 3D case.
    From top to bottom: Total solution time, final residual norm, and approximate contraction factor \(\hat{\rho}\).
    One can observe that for 64 and 512 subdomains, the asynchronous and the synchronous method take almost the same time for the solve, with a slight advantage for the asynchronous method.
    For 4096 subdomains, however, the synchronous method is seen to take significantly more time.
    The reason for this is not the solution of the coarse-grid problem, as in 2D, but the cost of the data exchange.
    The effect on the asynchronous method is much less pronounced.
  }
  \label{fig:weakScaling2level3d}
\end{figure}

The results of a weak scaling experiment are shown in \Cref{fig:weakScaling2level3d}.
We observe  behavior that is similar to the 2D case.
We notice however that the size of the coarse-grid problem (and hence the solution of the coarse-grid problem) are not the issue here.
At 4096 ranks, the coarse-grid problem is an order of magnitude smaller than the typical subdomain problem.
The apparent slowdown of the synchronous method is caused by the cost of exchanging information between the coarse grid and the subdomains.
While a slowdown is also visible in the asynchronous method, it is much less pronounced, resulting in a speedup of 4x over the synchronous method.

We also observe that compared to the 2D case where direct solvers were used, both synchronous and asynchronous iterations terminate almost exactly once the prescribed tolerance has been achieved.
The reason for this is that convergence checks occur much more frequently as the tolerance is reached, since the iterative sub-solves converge to their local tolerance typically within one iteration.

\section{Conclusion}

In the present work, we have explored the use of asynchronous alternatives to conventional (synchronous) one-level and two-level domain decomposition solvers.
To the best of our knowledge, we proposed the first truly asynchronous two-level method, where each processor can do different number of updates (iterations).
Several options to achieve asynchronous communication were tested, and we found that our use case benefited most from using \mpi{Win_lock_all} / \mpi{Win_unlock_all}, remote \mpi{Put}s and local \mpi{Get}s.
The numerical results presented demonstrate that asynchronous iterations can be considered a viable alternative to synchronous methods, despite partial availability of information from neighbors.
Asynchronous methods seem to be beneficial already at modest core count, even for load balanced scenarios.
In the presence of load imbalance, their performance advantage becomes even clearer, and we observed speedups up to 4x.
While we focused our attention on a particular Schwarz method, it is of inherent interest to explore asynchronous variants of other, potentially more effective domain decomposition methods involving deflation or non-overlapping decompositions (such as FETI and BDDC) or more than two levels.
The presented inclusion of a novel asynchronous coarse-grid correction paves the way for asynchronous methods to be used in extremely scalable parallel solvers.

\section*{Acknowledgment}

Sandia National Laboratories is a multimission laboratory managed and operated by National Technology and Engineering Solutions of Sandia, LLC, a wholly owned subsidiary of Honeywell International, Inc., for the U.S. Department of Energy’s National Nuclear Security Administration under contract DE-NA0003525.

This paper describes objective technical results and analysis. Any subjective views or opinions that might be expressed in the paper do not necessarily represent the views of the U.S. Department of Energy or the United States Government. \newline
SAND Number: SAND2020-8220 J

This material is based upon work supported by the U.S. Department of Energy, Office of Science, Office of Advanced Scientific Computing Research, Applied Mathematics program under Award Numbers DE-SC-0016564.
This research used resources of the National Energy Research Scientific Computing Center, a DOE Office of Science User Facility supported by the Office of Science of the U.S. Department of Energy under Contract No. DE-AC02-05CH11231.

\bibliographystyle{siamplain}
\bibliography{ref}

\end{document}